\documentclass{birkjour}
\usepackage{amsmath}
\usepackage{amssymb}
\newtheorem{thm}{Theorem}[section]
\newtheorem{cor}[thm]{Corollary}
\newtheorem{conj}[thm]{Conjecture}
\newtheorem{lem}[thm]{Lemma}
\newtheorem{prop}[thm]{Proposition}
\newtheorem{cons}[thm]{Construction}
\theoremstyle{definition}
\newtheorem{defn}[thm]{Definition}
\theoremstyle{remark}
\newtheorem{rem}[thm]{Remark}
\newtheorem{ex}[thm]{Example}
\newtheorem{exs}[thm]{Examples}
\long\def\Thm#1{\begin{thm} #1 \end{thm}}
\long\def\Cor#1{\begin{cor} #1 \end{cor}}

\long\def\Lem#1{\begin{lem} #1 \end{lem}}
\long\def\Prop#1{\begin{prop} #1 \end{prop}}

\long\def\Def#1{\begin{defn} #1 \end{defn}}
\long\def\Rem#1{\begin{rem} #1 \end{rem}}
\long\def\Ex#1{\begin{ex} #1 \end{ex}}
\long\def\Exs#1{\begin{exs} #1 \end{exs}}
\def\bar#1{\overline{#1}}
\def\Sect{\section}
\def\Rarr#1#2{\xrightarrow[#2]{#1}}

\long\def\Ref#1#2#3#4#5#6{
\bibitem{#1}
{\rm #2,}
\textit{#3.}
{\rm #4}
\textbf{#5}
{\rm #6.}
}
\long\def\Refb#1#2#3#4{
\bibitem{#1}
{\rm #2,}
\textit{#3.}
#4.
}
\def\Qq{{\mathbb Q}}%blackboard bold Q
\def\Zz{{\mathbb Z}}%blackboard bold Z
\def\Nn{{\mathbb N}}%blackboard bold Z
\def\Rr{{\mathbb R}}%blackboard bold R
\def\Cc{{\mathbb C}}%blackboard bold C
\def\Tt{{\mathbb T}}%blackboard bold T
\def\Ff{{\mathbb F}}%blackboard bold F

\def\Sq{{\rm Sq}}

\def\O{{\rm O}}
\def\U{{\rm U}}
\def\phi{\varphi}
\def\into{\hookrightarrow}
\def\iso{\cong}%isomorphism
\def\leq{\leqslant}
\def\geq{\geqslant}

\def\st{\mid}
\def\Zero{{\rm Zero}}

\def\cl#1{\bar{#1}}

\def\Sq{{\rm Sq}}
\begin{document}

\title{Some remarks on the parametrized Borsuk-Ulam theorem}

\author{Michael Crabb}
\address{%
Institute of Mathematics,\\
University of Aberdeen, \\
Aberdeen AB24 3UE, \\
UK}
\email{m.crabb@abdn.ac.uk}
\author{Mahender Singh}
\address{Indian Institute of Science Education and Research (IISER) Mohali,\\
Sector 81, Knowledge City,\\
SAS Nagar (Mohali), Post Office Manauli,\\
Punjab 140306,\\
India}
\email{mahender@iisermohali.ac.in}
\date{October 2017}
\begin{abstract}
Given a locally trivial fibre bundle $E\to B$ (with fibres and base finite
complexes), an orthogonal real line bundle $\lambda$ over $E$ 
and a real vector bundle $\xi$ over $B$, we consider a fibrewise map
$f: S(\lambda ) \to\xi$ over $B$ defined on the unit sphere bundle of 
$\lambda$. Following the fundamental work of Jaworowski
and Dold on the parametrized
Borsuk-Ulam theorem, we investigate lower bounds on the cohomological dimension
of the set $\{ v\in S(\lambda ) \st f(v)=f(-v)\}$. 
\end{abstract}
 %equivariant fibre spaces and bundles
\maketitle
\Sect{Introduction}
Let $B$ be a connected compact ENR (Euclidean Neighbourhood Retract), 
$\pi : E\to B$ a locally trivial fibre
bundle with compact ENR fibres, 
$\lambda$ a real line bundle over $E$ and $\xi$ a real vector 
bundle of dimension $n$ over $B$.
We may assume that $\lambda$ and $\xi$ are equipped with inner products and
write $S(\lambda )$ and $S(\xi )$ for their sphere bundles.

Suppose that $f : S(\lambda ) \to \xi$ is a fibrewise
map over $B$. 
We shall be concerned with the subset
$$
\tilde Z =\{ v\in S(\lambda )\st f(v)=f(-v)\}
$$
of $S(\lambda )$ and its image $Z$ in $E$ under the projection $S(\lambda )
\to E$.
In particular, if $f(-v)=-f(v)$ for all $v\in S(\lambda)$, then
$\tilde Z$ is the zero-set $\{ v\in S(\lambda ) \st f(v)=0\}$ of $f$.

Throughout the paper we shall use representable cohomology
$H^*$ (as in \cite[Section 8]{borsuk}), usually with $\Ff_2$-coefficients.
Our goal is to estimate the size of the space $\tilde Z$ by giving
a lower bound on its cohomological dimension.
Early estimates of this type were obtained by Jaworowski \cite{Jan} and
Dold \cite{Dold}. 
We present here two theorems that generalize the result
\cite[Corollary 1.5]{Dold} in Dold's fundamental paper on the parametrized
Borsuk-Ulam theorem.

The first is quite elementary.
But it leads to easy proofs of a number of results in the recent literature
(\cite[Corollary 1.5]{KM}, \cite[Theorem 1.3]{MPSS},
\cite[Theorem 1.5]{MPS}, \cite[Theorem 1.3]{S2}). 
In the statement below, $e(\lambda )=w_1(\lambda )
\in H^1(E;\,\Ff_2)$ is the $\Ff_2$-Euler class of $\lambda$.
\Thm{\label{stiefel}
Suppose that $H^i(B;\,\Ff_2)=0$ for $i>d$ and that, for some $k\geq n$,
there is a class $b\in H^d(B;\,\Ff_2)$ such that 
$\pi^*(b)\cdot e(\lambda )^k\not=0$. 
Then the restriction of $\pi^*(b)\cdot e(\lambda )^{k-n}
\in H^{d+k-n}(E;\,\Ff_2)$ to $Z$ is non-zero.

It follows that $H^j(\tilde Z;\,\Ff_2)$ is non-zero for some
$j\geq d+k-n$.
}
The second is rather deeper and leads to stronger results in specific examples.
Given a real vector bundle $\eta$, we shall write $P(\eta )$ for its
projective bundle and $H$ for the Hopf line bundle over $P(\eta )$.
\Thm{\label{mainthm}
Suppose that, for some $r\geq 1$, $\pi : E\to B$ admits a factorization
$$
E=S(\zeta_r)\to S(\zeta_{r-1})\to \cdots \to S(\zeta_1)\to P(\eta )\to B,
$$
where $\eta$ is an $(m+1)$-dimensional real vector bundle over $B$,
and $\zeta_i$, for $i=1,\ldots ,r$, is a real vector bundle of
dimension $l_i+1$ over $P(\eta )$ if $i=1$, 
over $S(\zeta_{i-1})$ if $i>1$,
with $\Ff_2$-Euler class $e(\zeta_i)=w_{l_i+1}(\zeta_i)$ equal to zero.
The fibres of $\pi$ are thus manifolds of dimension $l+m$,
where $l=l_1+ \ldots + l_r$.
Suppose further that $\lambda$ is the pullback of $H$ over $P(\eta )$.

Let $d$ be maximal such that $H^d(B;\,\Ff_2)\not=0$. 
Then, if $n\leq m$, the cohomology group
$H^{d+l+m-n}(Z;\,\Ff_2)$ is non-zero, and hence  
$H^j(\tilde Z;\,\Ff_2)$ is non-zero for some $j\geq d+l+m-n$.
}
Two cases in which these conditions are satisfied are described
in Propositions \ref{ms} and  \ref{mpss}, which pursue ideas introduced
in \cite{S2} and \cite{MPSS}.

The proofs of Theorems \ref{stiefel} and \ref{mainthm}
are given in Sections 2 and 3. An analogous theory for complex vector
bundles is sketched in Section 4.
In Section 5 we discuss the extension of the parametrized Borsuk-Ulam 
theorem for a sphere bundle $S(\xi )$ with the antipodal involution to the
case of a spherical fibration with a fibre-preserving free involution.
Necessary material on the Euler class of a spherical fibration is included
as an appendix (Section 6).
\Sect{An elementary condition}
The basic Borsuk-Ulam Theorem 
as formulated, for example, in \cite[Proposition 2.7]{borsuk}
specializes to the following proposition, which will be fundamental
to our discussion.
To be precise we should write $\lambda\otimes\pi^*\xi$ in its statement,
rather than $\lambda\otimes\xi$.
In order to simplify notation we shall often, as here,
use the same symbol for a bundle and its pullback to some space.
\Prop{\label{borsuk}
Suppose that $a\in H^i(E;\,\Ff_2 )$ is a cohomology class 
such that $a\cdot e(\lambda\otimes\xi )\not=0\in H^{i+n}(E;\, \Ff_2)$.
Then $a$ restricts
to  a non-zero class in $H^i(Z ;\,\Ff_2 )$. 
}
\begin{proof}
The map $f$ determines a section $s$ of $\lambda\otimes \xi$ over $E$:
$$
s(x) = v\otimes (f(v)-f(-v)) \text{\ for\ } v\in S(\lambda_x)
$$
with zero-set $\{ x\in E \st s(x)=0\}$ equal to the subset $Z$. 
The conclusion follows at once from \cite[Proposition 2.7]{borsuk}
or its extension to spherical fibrations given as Proposition \ref{loc}
in the Appendix. 
\end{proof}
\Cor{Under the conditions of Proposition \ref{borsuk}, 
$H^j(\tilde Z;\,\Ff_2)$ is non-zero for some $j\geq i$.
}
\begin{proof}
Notice that $\tilde Z =S(\lambda\, |\, Z)$.
There is, thus, a long exact sequence
$$
\cdots\to H^j(Z;\,\Ff_2) \to H^j(\tilde Z;\,\Ff_2) \to
H^j(Z;\,\Ff_2) \Rarr{e(\lambda )\cdot}{} H^{j+1}(Z;\,\Ff_2)\to\cdots
$$
The result follows easily from the fact that $e(\lambda )$ is nilpotent.
\end{proof}
\begin{proof}[Proof of Theorem \ref{stiefel}]

Recall that the mod $2$ Euler class $e(\lambda\otimes\xi )$ is given by
$$
w_n(\lambda \otimes\xi ) =\sum_{i=0}^n e(\lambda )^i w_{n-i}(\xi ).
$$
Because $b\cdot w_{n-i}(\xi )=0$ for $i<n$,
we have
$\pi^*(b)\cdot e(\lambda )^{k-n}\cdot e(\lambda\otimes\xi )
=\sum_{i=0}^n \pi^*(b)\cdot e(\lambda )^{k-n}\cdot e(\lambda )^iw_{n-i}(\xi )
= \pi^*(b)\cdot e(\lambda )^k$, which is non-zero. 
Take $a=\pi^*(b)\cdot e(\lambda )^{k-n}$ in Proposition \ref{borsuk}.
\end{proof}
\Rem{If we replace the condition $\pi^*(b)\cdot e(\lambda )^k\not=0$ by
$\pi^*(b)\cdot c\cdot e(\lambda )^n\not=0$ for some class 
$c\in H^{k-n}(E;\,\Ff_2)$, the same argument shows that $\pi^*(b)\cdot c$
restricts to a non-zero class in $H^{d+k-n}(Z;\,\Ff_2)$.
}
\Cor{\label{simple}
Suppose that $d\geq 0$ is maximal such that $H^d(B;\,\Ff_2)$ is non-zero
and that, for any fibre $F$ of $\pi$,
the restriction $H^*(E;\,\Ff_2)\to H^*(F;\,\Ff_2)$ is surjective and 
$e(\lambda | F)^k\in H^k(F;\,\Ff_2)$ is non-zero.
Then $H^{d+k-n}(Z;\,\Ff_2)$ is non-zero.
}
\begin{proof}
In this case multiplication by $e(\lambda )^k : H^*(B) \to H^*(E)$
is injective, by the Leray-Hirsch Theorem, which says that the bundle
$E\to B$ is `homologically trivial'.
\end{proof}
\Exs{\label{exs}
We give some examples from \cite{S2, KM, MPSS, MPS}.

\par\noindent (i). (\cite[Theorem 1.3]{S2}).
Suppose that the restriction $S(\lambda |F) \to F$ of the bundle $S(\lambda )
\to E$ to a fibre $F$ is homeomorphic to the projection 
$$
V_r(\Rr^{r+s})\to V_r(\Rr^{r+s})/\{\pm 1\}
$$
from the Stiefel manifold of orthogonal $r$-frames in $\Rr^{r+s}$
to the projective Stiefel manifold.
Then we may take $k$ to be the smallest integer in the range $s\leq k <r+s$ 
such that $\binom{r+s}{k+1}$ is odd. (See 
%\cite[Theorem 3.1]{S2},
\cite[Theorem 1.6]{GH}.)

\par\noindent (ii). (\cite[Corollary 1.5]{KM}, \cite[Theorem 1.3]{MPSS}).
Suppose that $S(\lambda |F )\to F$ is homeomorphic to the projection
$$
S^{n_1}\times\cdots \times S^{n_r} \to (S^{n_1}\times\cdots\times S^{n_r})/
\{\pm 1\}
$$ 
from a product of $r$ spheres to the quotient by the involution
that acts antipodally on each factor.
We may take $k$ to be $\min\{ n_1,\ldots ,n_r\}$.
(See 
%\cite[Theorem 3.1]{MPSS}, 
\cite[Theorem 2.1]{DMD}.)

\par\noindent (iii). (\cite[Theorem 1.5]{MPS}).
Let $r,\, s > 1$ be integers.
Suppose that the fibre $S(\lambda | F)$ has cohomology ring
$H^*(S(\lambda | F);\,\Ff_2)
=\Ff_21\oplus\Ff_2u\oplus \Ff_2v\oplus \Ff_2w$,
where $u,\, v,\, w$ have degrees $r$, $r+s$, $2r+s$, respectively
and $uv=0$.
Then we may take $k=2r+s$.
(See 
%\cite[Theorem 2.4]{MPS}, 
\cite[Theorem 4.1]{PSS}, at least for the case $r=s$.)
}
We shall return to the examples (i) and (ii) in the next section.
\Rem{Here is an illustration of the example (iii).
Let $\alpha$, $\beta$ and $\gamma$ be real vector bundles over $B$
of dimension $r$, $s+1$ and $r+s+1$, respectively. We form the real projective
bundles $P(\alpha \oplus\beta )\to B$ and $P(\alpha\oplus\gamma )\to B$.
Both contain $P(\alpha )$ as a subbundle.
We define $E\to B$ by gluing the two bundles along $P(\alpha )$:
$$
E = P(\alpha \oplus\beta )\cup_{P(\alpha )} P(\alpha\oplus\gamma ) \to B
$$
and take $\lambda$ to be the line bundle that restricts to the Hopf line
bundle on each of the subspaces $P(\alpha\oplus\beta )$ and 
$P(\alpha\oplus \gamma )$.
The double cover $S(\lambda )$ may be described by a gluing construction
or as a fibrewise join:
$$
S(\lambda ) = S(\alpha\oplus\beta )\cup_{S(\alpha )} S(\alpha\oplus\gamma )
=S(\alpha )*_B (S(\beta ) \sqcup S(\gamma )).
\qed
$$
}
To be accurate,
the papers \cite{S2, KM, MPSS, MPS} work in a rather more general setting
than that considered here: the base $B$ is not required to be 
a compact ENR and the involution on the vector bundle $\xi$ is not 
necessarily antipodal.
We shall not deal with more general base spaces here, but we explain in
Section 5 how to incorporate more general involutions into the theory.

For the covering dimension there is a simpler estimate.
\Cor{Suppose that $B$ is a closed smooth manifold of dimension $d$ and
that, for some point $b\in B$, the restriction
$e(\lambda |F)^k\in H^k(F;\,\Ff_2)$ to the fibre $F$ of $\pi$ at $b$
is non-zero.
Then the covering dimension of $Z$ and of $\tilde Z$ is at least $d+k-n$.
}
\begin{proof}
We choose an embedding $D(\Rr^d)\into B$ of a closed $d$-disc centred at $b$.
The restriction of the bundle $E\to B$ to $D(\Rr^d)$ is trivial. 
Now we map the $d$-dimensional sphere $B'=S(\Rr\oplus \Rr^d)$ to
$D(\Rr^d)$, and so to $B$, by projection $(t,v)\mapsto v$.
The bundles $E\to B$, $\lambda$ over $E$ and $\xi$ over $B$ pullback to 
bundles $E'\to B'$, $\lambda$ over $E'$ and $\xi'$ over $B'$, and the
map $f$ lifts to $f' : S(\lambda')\to \xi'$ with the associated
subset $Z'\subseteq E'$ lifting $Z\subseteq E$.

We can apply Corollary \ref{simple} to $f'$ to deduce that 
$H^{d+k-n}(Z';\,\Ff_2)$ is non-zero and so that the covering dimension
of $Z'$ is at least $d+k-n$.
But $Z'$ is the union of two closed subspaces each homeomorphic to 
$Z|D(\Rr^d)$. Hence $Z|D(\Rr^d)$ and so $Z$ have covering dimension
greater than or equal to $d+k-n$.
\end{proof}
A more general result on the covering dimension is described in 
Proposition \ref{cover}.
\Sect{Sphere bundles}
Consider an $m+1$-dimensional real vector bundle $\eta$ over $B$ 
with projective bundle $P=P(\eta ) \to B$.

It is convenient to introduce the polynomials
$$
p(T)=T^{n}+w_1(\xi )T^{n-1} + \ldots + w_{n}(\xi )
$$
and
$$
q(T)=T^{m+1}+w_1(\eta )T^{m} + \ldots + w_{m+1}(\eta )
$$
in $H^*(B;\,\Ff_2)[T]$.
\Rem{\label{deriv}
Recall that 
$$
p(T+w_1(\lambda)) =\sum_{i=0}^{n} w_{n-i}(\lambda\otimes\xi ) T^i
\in H^*(E;\, \Ff_2)[T].
$$
}
\Rem{\label{div}
The Euler class $e(H\otimes\xi )\in H^n(P;\,\Ff_2)$ (of the tensor product of
the Hopf line bundle $H$ and the pullback of $\xi$) is zero if
and only if $q(T)$ divides $p(T)$ in $H^*(B;\,\Ff_2)$.
}
Suppose that $\rho :E \to P$ is a fibre bundle with structure
group a compact Lie group $G$ acting smoothly on a closed connected
manifold $M$ of dimension $l$. Thus $E = Q\times_G M$ for a principal
$G$-bundle $Q\to P$. (More generally, we could just take $E \to P$
to be a fibrewise manifold with fibre $M$.)
We take $\pi$ to be the composition $E\to P\to B$ and
$\lambda =\rho^*H$ to be the pullback of the Hopf line bundle.
\Lem{\label{split}
Suppose that there is a class $\sigma\in H^l(E;\,\Ff_2)$ restricting
to the generator of the cohomology $H^l(M;\,\Ff_2)=\Ff_2$ of a fibre.
Then, for each $i\geq 0$, the homomorphism
$$
x\mapsto \rho^*(x)\cdot\sigma : H^i(P;\,\Ff_2) \to H^{i+l}(E;\,\Ff_2)
$$
is a split injection.
}
\begin{proof}
We have a fibrewise Umkehr map $\rho_! : H^{i+l}(E;\,\Ff_2) \to H^i(P;\,\Ff_2)$
and $\rho_!(\rho^*(x)\cdot \sigma ) = x\cdot\rho_!(\sigma )=x$ for
$x\in H^i(P;\,\Ff_2)$, because $\rho_!(\sigma )=1\in H^0(P;\,\Ff_2)$.
\end{proof}
\Rem{The class $\sigma$ should be thought of as a `homology section'
of $\rho : E \to P$. A genuine section $s : P \to E$
determines a class $\sigma = s_!(1)\in H^l(E;\,\Ff_2)$.
}
\Prop{\label{main}
Suppose that there is a class $\sigma\in H^l(E;\,\Ff_2)$ satisfying the
condition of Lemma \ref{split} and that $n\leq m$.
Let $b\in H^d(B;\,\Ff_2)$ be a non-zero class of maximal degree.
Then $\pi^*(b)\cdot e(\lambda )^{m-n}\cdot \sigma
\in H^{d+m-n+l}(E;\,\Ff_2)$ restricts to a non-zero class
in $H^{d+m-n+l}(Z;\,\Ff_2)$.
}
\begin{proof}
Write $a=\pi^*(b)\cdot e(\lambda )^{m-n}\cdot \sigma$.
Then
$$
a\cdot e(\lambda\otimes\xi ) = 
\rho^*(b\cdot e(H)^{m-n}\cdot e(H\otimes\xi ))\cdot\sigma .
$$
But $b\cdot e(H)^{m-n}\cdot e(H\otimes\xi )
\in H^{d+m}(P;\,\Ff_2)$ is non-zero.
So the result follows from Proposition \ref{borsuk}.
\end{proof}
\Rem{When $\rho$ is the identity $E=P\to P$, this reduces to the original result
of Dold \cite[(1.7)]{Dold}.
}
Let $\zeta$ be a real vector bundle of dimension $l+1$ over $P$.
We consider the bundle $E=S(\zeta )\to P$.
\Lem{\label{gen}
There is a class $\sigma \in H^l(S(\zeta );\,\Ff_2)$
restricting to the generator of the cohomology $H^l(M;\,\Ff_2)$ of
of a fibre if and only if the Euler class $e(\zeta )\in H^{l+1}(P;\,\Ff_2)$
is zero.
}
\begin{proof}
We have an exact Gysin sequence
of the pair $(D(\zeta), S(\zeta))$:
$$
H^{-1}(P;\,\Ff_2)=0 \to H^{l}(P;\,\Ff_2) \to 
H^{l}(S(\zeta);\,\Ff_2)
$$
$$
\to H^0(P;\,\Ff_2)=\Ff_2 \Rarr{e(\zeta )}{} H^{l+1}(P;\,\Ff_2),
$$
where we have used the Thom isomorphism 
$$
H^{l+1+*}(D(\zeta),S(\zeta);\,\Ff_2) \iso H^*(P;\,\Ff_2),
$$
so that the restriction $H^*(D(\zeta),S(\zeta);\,\Ff_2)
\to H^*(D(\zeta))$ corresponds to multiplication by the
Euler class $e(\zeta)=w_{l+1}(\zeta )$: 
$H^*(P;\,\Ff_2) \to H^{*+l+1}(P;\,\Ff_2)$.

Thus there is a lift $\sigma \in H^{l}(S(\zeta);\,\Ff_2)$
of $1\in H^0(P;\,\Ff_2)$ if and only if $e(\zeta )=0$.
\end{proof}
\Rem{\label{sq}
The class $\sigma$ satisfies 
$\Sq (\sigma) - w(\zeta )\sigma \in H^*(P;\,\Ff_2)$, 
because the total Steenrod square $\Sq$ acts on the Thom class of $\zeta$ as
multiplication by the total Stiefel-Whitney class.
In particular, we have the identity 
$$
\sigma^2 - w_l(\zeta )\sigma \in H^{2l}(P;\,\Ff_2).
$$
}
\Ex{\label{two}
(\cite{S2}).
Consider the bundle of Stiefel manifolds $\tilde E = \O (\Rr^2,\eta )$,
with fibre at $x\in B$ the space $V_2(\eta_x)=\O (\Rr^2,\eta_x)$
of orthogonal $2$-frames in $\eta_x$, 
equipped with the involution $-1$ on
$\eta$. The quotient $E$ is a bundle of projective Stiefel manifolds. 
Restriction to the first factor $\Rr$ in $\Rr^2$ gives a map
$\O (\Rr^2 ,\eta )\to \O (\Rr ,\eta )=S(\eta )$.
Then we can express $E$ as the sphere bundle $S(\zeta )\to
P(\eta )$, where $\zeta$ is the tensor product $H\otimes H^\perp$
of $H$ and its orthogonal complement $H^\perp$ in the
pullback of $\eta$. Thus $l=m-1$.
The line bundle $\lambda$ is the pullback of $H$.

Now 
$$
w_{l+1}(\zeta )
=(m+1)w_1(H)^m + mw_1(\eta )w_1(H)^{m-1} + \ldots + w_m(\eta )
$$
by the formula $w_m(H\otimes\eta )=q'(w_1(H))$.
(See Remark \ref{deriv}.)
For $w_{l+1}(\zeta )$ to vanish it is necessary and sufficient
that $m+1$ be even and $w_i(\eta )=0$ for $i$ odd.
This is true if $\eta$ admits a complex structure,
and in that special case $\zeta$ has a trivial $1$-dimensional summand.
}
\Ex{Suppose that $\eta$ has a complex structure, so that $m+1$ is even.
Consider the complex Stiefel bundle $\U (\Cc^2,\eta )$ with the involution
$-1$ and quotient $E$.
We can express $E$ as $S(\zeta )\to P=P(\eta )$, where $\zeta =H\otimes
(\Cc\otimes H)^\perp$, so that $l=m-2$.
By Remark \ref{deriv} again, we have
$$
\textstyle
w_{l+1}(\zeta )
=\binom{m+1}{2}w_1(H)^{m-1} + \binom{m}{2}w_1(\eta )w_1(H)^{m-2} + \ldots + 
\binom{2}{2}w_{m-1}(\eta ),
$$
which vanishes if and only if $m+1$ is divisible by $4$ and
$w_{2i}(\eta )=0$ for $i$ odd.
This holds if $\eta$ admits a quaternionic structure.
}
Consider, more generally, $r$ vector bundles $\zeta_i$ of dimension $l_i+1$,
$i=1,\ldots ,r$, over $P$. 
Take $E = S(\zeta_1)\times_P \cdots\times_P S(\zeta_r)$ and
let $\lambda$ be the pullback of $H$ over $P$.
Write $l=l_1+\ldots +l_r$.
\Lem{\label{euler}
There is a class $\sigma \in H^l(E;\,\Ff_2)$ restricting to
a generator of the cohomology $H^l(M;\,\Ff_2)$ of a fibre if and only
if $e(\zeta_i)=0\in H^{l_i+1}(P;\,\Ff_2)$ for $i=1,\ldots ,r$.
}
\begin{proof}
If each $e(\zeta_i)$ is zero, choose $\sigma_i\in H^{l_i}(S(\zeta_i);\,\Ff_2)$
as in Lemma \ref{gen}. Then take $\sigma =\sigma_1\cdot \ldots\cdot\sigma_r$.

Conversely, given $\sigma$, we can produce a class $\sigma_i$ as the image
of $\sigma$ under the Umkehr homomorphism
$$
(\pi_i)_! : H^l(S(\zeta_1)\times_P \cdots\times_P S(\zeta_r);\,\Ff_2)
\to H^{l_i}(S(\zeta_i);\,\Ff_2)
$$
of the fibrewise projection 
$\pi_i : S(\zeta_1)\times_P \cdots\times_P S(\zeta_r) \to S(\zeta_i)$ over $P$.
\end{proof}
A special case of our main result Theorem \ref{mainthm}
follows immediately from Proposition \ref{main}.
\Prop{\label{mainprop}
Suppose that $e(\zeta_i)=0$ for $i=1,\ldots ,r$, that 
$d$ is maximal such that $H^d(B;\,\Ff_2)\not=0$ and that $n\leq m$.
Then  $H^{d+m-n+l}(Z;\,\Ff_2)$ is non-zero.
\qed
}
\Rem{By Leray-Hirsch, $H^*(E;\,\Ff_2)$ is a free $H^*(P;\,\Ff_2)$-module 
on the classes
$\sigma_{i_1}\cdots \sigma_{i_k}$, $1\leq i_1< \ldots < i_k\leq r$,
where $0\leq k\leq r$.
Remark \ref{sq} allows us to describe $H^*(E;\,\Ff_2)$ as
$$
H^*(P;\,\Ff_2)[V_1,\ldots ,V_r]/(V_i^2-w_{l_i}(\zeta_i)V_i-s_i\, |\, 
i=1,\ldots ,r),
$$
where $V_i$ has dimension $l_i$
and $s_i=\sigma_i^2 -w_{l_i}(\zeta_i)\sigma_i\in H^{2l_i}(P;\,\Ff_2)$.
(Compare \cite[Section 4]{MPSS}.)

This decomposition is a homological version of the stable splitting in the next
Remark \ref{davis}.
}
\Rem{\label{davis}
The condition of Lemma \ref{euler} is satisfied if each $S(\zeta_i)$
has a section over $P$ and we can then split
$\zeta_i$ as $\Rr\oplus\nu_i$ and identify $S(\zeta_i)$ with
the fibrewise one-point compactification $(\nu_i)_P^+$ over $P$.
Now the fibrewise suspension $\Sigma_P (S(\zeta)_{+P})$ 
splits as a fibrewise
wedge $\Sigma_P(P\times S^0)\vee_P\Sigma_P(\nu_i)^+_P$.
It follows that we have a fibrewise homotopy equivalence
$$
\Sigma_P(E_{+P}) \simeq
\Sigma_P\left( ((P\times S^0) \vee_P (\nu_1)^+_P)\wedge_P \cdots \wedge_P
((P\times S^0) \vee_P (\nu_r)^+_P)\right).
$$
This then permits us,
using the notation $P^\eta$ for the Thom space of a vector bundle
$\eta$ over $P$, to describe $\Sigma (E_+)$ as a wedge of 
$2^{r}$ suspensions of Thom spaces over $P$:
$$
\Sigma (E_+) \simeq
\bigvee_{1\leq i_1 < i_2 < \ldots < i_k\leq r} 
\Sigma (P^{\nu_{i_1}\oplus \cdots\oplus \nu_{i_k}})
$$
where $0\leq k \leq r$.
There is thus a homotopy equivalence
$$
\Sigma^r E_+ \simeq
\bigvee_{1\leq i_1 < i_2 < \ldots < i_k\leq r} 
\Sigma^{r-k}(P^{\zeta_{i_1}\oplus \cdots\oplus \zeta_{i_k}})\, .
$$
}
\begin{proof}[Proof of Theorem \ref{mainthm}]
Choose $\sigma_i \in H^{l_i}(S(\zeta_i);\,\Ff_2)$, for $i=1,\ldots ,r$, 
as in Lemma \ref{gen}, and define $\sigma \in H^l(E;\,\Ff_2)$ to be
the product of their lifts to $E$. 
Then, as in Lemma \ref{split},
$x\mapsto \rho^*(x)\cdot\sigma : H^i(P;\, \Ff_2) \to H^{i+l}(E;\,\Ff_2)$ 
is a composition of split injections
$$
H^i(P;\Ff_2) \to H^{i+l_1}(S(\zeta_1);\Ff_2) \to \cdots
\to H^{i+l_1+\ldots  +l_r}(S(\zeta_r);\,\Ff_2).
$$
The assertion of the theorem follows from Proposition \ref{main}.
\end{proof}
This allows us to extend Example \ref{two} to higher dimensions.
\Prop{\label{ms}
{\rm (See \cite{S2}.)}
Fix an integer $r\geq 1$ and let 
$s\in\Nn$ be defined by $2^{s-1}\leq r <2^s$. 

Let $\eta$ be a real vector bundle of dimension $m+1\geq r+1$ over $B$.
Consider the bundle of
Stiefel manifolds $\tilde E=\O (\Rr^{r+1},\eta )$ with the free involution
$-1$.
Let $\pi : E\to B$ be the quotient by the action of $\O (1)$
and let $\lambda$ be the associated real line bundle.

Suppose that 
$m+1$ is divisible by $2^s$ and, for $j>0$,
$w_j(\eta )=0$ if $j$ is not divisible by $2^s$.

Then $\pi : E\to B$ and $\lambda$ satisfy the conditions of 
Theorem \ref{mainthm} with $l_i=m-i$ and $l+m = (r+1)m-r(r+1)/2$.
}
\begin{proof}
Restriction to the first factor gives a map $\O (\Rr^{r+1},\eta )\to
S(\eta )$ and so, on quotients, $E\to P(\eta )$. 
We can then interpret $E\to P(\eta )$
as the Stiefel bundle $\O (\Rr^r,H\otimes H^\perp) \to P(\eta )$,
where $H^\perp$ is the orthogonal complement of $H$ in $\eta$.
With $\zeta_1=H\otimes H^\perp$, we take $\zeta_i$ to be the tangent bundle
of $S(\zeta_{i-1})$ for $i=2,\ldots ,r$. Thus, $\Rr\oplus \zeta_i=\zeta_{i-1}$,
$\Rr^{i-1}\oplus\zeta_i=H\otimes H^\perp$ and 
$\Rr^i\oplus\zeta_i = H\otimes \eta$.
$\zeta_i$ of dimension $m+1-i$ and $e(\zeta_i)=w_{m+1-i}(\zeta_i)=
w_{m+1-i}(H\otimes\eta )$.

We have $e(\zeta_i) =w_{m-i+1}(H\otimes\eta )$.
By Remark \ref{deriv},
the condition that $e(\zeta_i)$ should vanish for $i=1,\ldots ,r$ is that
$\binom{m+1-j}{i} w_j(\eta )=0$ for $1\leq i\leq r$ and $0\leq j\leq m+1-i$. 
This holds if and only if $m+1$ is divisible by $2^s$ and $w_j(\eta )=0$
if $j$ is not divisible by $2^s$.
%Write $m+1-j = 2^t({\rm odd})$. 
\end{proof}
Let $\mu$ be a real vector bundle of dimension $l+1$ over $B$.
We now consider $\zeta =H\otimes\mu$ and write
$$
r(T) =T^{l+1}+w_1(\mu )T^l + \ldots + w_{l+1}(\mu )\in H^*(B;\,\Ff_2)[T]
$$
\Lem{\label{spherebundle}
There is a class $\sigma \in H^l(S(H\otimes\mu );\,\Ff_2)$
restricting to the generator of the cohomology $H^l(M;\,\Ff_2)$ of
of a fibre if and only if the polynomial $r(T)$ is divisible by  $q(T)$.
In that case, $\sigma^2 - r'(w_1(H))\sigma \in H^{2l}(P;\,\Ff_2)$.
}
Notice that this condition requires the dimensional restriction $l\geq m$ and
is satisfied if $\mu$ admits $\eta$ as a subbundle.
\begin{proof}
This is contained in Remark \ref{div}.
\end{proof}
\Prop{\label{mpss}
{\rm (See \cite{KM, MPSS}.)}
Let $\eta$ be a real vector bundle over $B$ of dimension $m+1$
and $\mu_1,\ldots ,\, \mu_r$, $r\geq 1$, real vector
bundles over $B$ with $\dim \mu_i =l_i+1\geq m+1$.
Consider the bundle $\tilde E = S(\eta )\times_B S(\mu_1)\times_B 
\cdots\times_B S(\mu_r)$ with the free diagonal antipodal involution.
Let $\pi : E\to B$ be the quotient by the action of $\O (1)$ and
let $\lambda$ be the associated real line bundle, 
so that $\tilde E=S(\lambda )$.

Suppose that
$$
r_i(T)=T^{l_i+1}+w_1(\mu_i)T^{l_i-1} + \ldots + w_{l_i+1}(\mu_i)\in 
H^*(B;\,\Ff_2)[T]
$$
is divisible by 
$q(T)=T^{m+1}+w_1(\eta )T^m + \ldots + w_{m+1}(\eta )$ 
for each $i=1,\ldots ,r$.

\smallskip

Then $\pi : E\to B$ and $\lambda$ satisfy the conditions of
Theorem \ref{mainthm}.
}
\begin{proof}
The bundle $E \to B$ can be identified with
$$
E=S(H\otimes\mu_1)\times_P \cdots \times_P S(H\otimes\mu_r) \to P \to B
$$
and $\lambda$ with the pullback of the Hopf bundle $H$.
\end{proof}
Thus Theorem \ref{mainthm}, or Proposition \ref{mainprop},
extends \cite[Theorem 1.4]{MPSS} to bundles.
\Rem{The condition of Proposition \ref{mpss}
is satisfied if each $\mu_i$ admits $\eta$ as a subbundle,
or equivalently if the sphere bundle $S(\zeta_i)$ of $\zeta_i=H\otimes\mu_i$
admits a section.
In that case, Remark \ref{davis} 
allows us to describe $\Sigma (E_+)$ as a wedge of 
$2^{r}$ Thom spaces over $P$:
$$
\Sigma (E_+) \simeq
\bigvee_{1\leq i_1 < i_2 < \ldots < i_k\leq r} 
\Sigma P^{\nu_{i_1}\oplus \cdots\oplus \nu_{i_k}}
$$
where $0\leq k \leq r$.
There is thus a homotopy equivalence
$$
\Sigma^r E_+ \simeq
\bigvee_{1\leq i_1 < i_2 < \ldots < i_k\leq r} 
\Sigma^{r-k}P^{H\otimes (\mu_{i_1}\oplus \cdots\oplus \mu_{i_k})}\, .
$$
(Compare the argument of Davis in \cite{DMD}.)
}
\Sect{Complex versions}
The methods extend readily to the complex theory.
We suppose now that $\lambda$ and $\xi$ are complex vector bundles
over $B$ of dimension $1$ and $n$.
From a map $f : S(\lambda ) \to \xi$ we construct a section $s$ of the
complex tensor product $\lambda^*\otimes\xi$:
$$
s(x) = v^* \otimes \int_\Tt zf(z^{-1}v)\text{\ for\ } v\in S(\lambda_x),
$$
where $v^*\in S(\lambda_x^*)$ is the dual generator and the integral
is over the circle group $\Tt =\{ z\in\Cc \st |z|=1\}$ (with Haar measure). 
We write $Z=\Zero (s)$, and then $\tilde Z = S(\lambda | Z)$ is
the set of points $v$ where $\int zf(z^{-1}v)$ is zero.

We use cohomology with $\Zz$-coefficients. The Euler classes are
$e(\lambda )= c_1(\lambda )$ and $e(\lambda^*\otimes\xi )
=\sum (-1)^ie(\lambda )^i c_{n-i}(\xi )$.
\Thm{\label{cxstiefel}
Suppose that $H^i(B;\,\Zz )=0$ for $i>d$ and that, for some $k\geq n$,
there is a class $b\in H^d(B;\,\Zz )$ such that 
$\pi^*(b)\cdot e(\lambda )^k\not=0$.
Then the restriction of $\pi^*(b)\cdot e(\lambda )^{2k-2n}
\in H^{d+2k-2n}(E;\,\Zz )$ to $Z$ is non-zero.

It follows that $H^j(\tilde Z;\,\Zz )$ is non-zero for some
$j\geq 1+d+2(k-n)$.
\qed
}
There are corresponding results for cohomology with coefficients
in $\Ff_p$ ($p$ a prime) or $\Qq$.
We state a corollary for $\Ff_p$-cohomology.
\Cor{Suppose that $d\geq 0$ is maximal such that $H^d(B;\,\Ff_p)$
is non-zero and that, for any fibre $F$ of $\pi$, the restriction
$H^*(E;\,\Ff_p)\to H^*(F;\,\Ff_p)$ is surjective and 
$e(\lambda |F)^k\in H^{2k}(F;\,\Ff_p)$ is non-zero. Then
$H^{d+2k-2n}(Z;\,\Ff_p)$ is non-zero, and
$H^j(\tilde Z;\,\Ff_p)$ is non-zero for some $j\geq 1+d+2(k-n)$.
}
\Ex{(\cite[Theorem 1.4]{S2}).
Suppose that the restriction $S(\lambda | F) \to F$ of the bundle 
$S(\lambda )\to E$ to a fibre $F$ is homeomorphic to the projection
$$
V^\Cc_r(\Cc^{r+s}) \to V^\Cc_r(\Cc^{r+s})/\Tt
$$
from the complex Stiefel manifold of unitary $r$-frames in $\Cc^{r+s}$ to the
projective complex Stiefel manifold.
We can take $k$ to be the least integer in the range $s\leq k<r+s$ such that
$\binom{r+s}{k+1}$ is not divisible by $p$.
(See 
%\cite[Theorem 3.2]{S2}
\cite[Theorems 1.1 and 1.2]{astey}.)
}
\Thm{\label{cxmainthm}
Suppose that $\eta$ is a complex vector bundle of dimension $m+1$ over
$B$
and that $\zeta_i$, for $i=1,\ldots ,r$, is a complex vector bundle
of dimension $l_i+1$ over $\Cc P(\eta )$ for $i=1$ and
over $S(\zeta_{i-1})$ for $i>1$, with each $\Zz$-Euler class
$e(\zeta_i)=c_{l_i+1}(\zeta_i)$ zero.
Let $E\to B$ be the bundle
$$
E = S(\zeta_r)\to S(\zeta_{r-1})\to \cdots \to S(\zeta_1) \to \Cc P(\eta ) \to B
$$
with $\lambda$ the pullback of the complex Hopf bundle $H$ over $\Cc P(\eta )$.

Let $d$ be maximal such that $H^d(B;\,\Zz )\not=0$.
Then, if $n\leq m$, the group $H^{d+2m-2n+2l}(Z;\,\Zz )$ is non-zero.

It follows that $H^j(\tilde Z;\,\Zz )$ is non-zero for some
$j\geq 1+d+2(m-n+l)$.
\qed
}
\Ex{Let $E=\U (\Cc^2,\eta )$, where $\eta$ is a complex bundle of even
dimension $m+1$ such that $c_i(\eta )=0$ for all odd $i$.
We take $\zeta$ to be $H^*\otimes H^\perp$ over $\Cc P(\eta )$,
where $H^\perp$ is the orthogonal complement of $H\subseteq \eta$, so that
$l=m-1$.
The condition on Chern classes is satisfied if the complex structure
on $\eta$ extends to a quaternionic structure.
}
\Sect{Spherical fibrations}
We sketch a generalization of the theory in which the sphere bundle
$S(\xi )$ with the antipodal involution is replaced by a fibrewise 
$\O (1)$-space $\Xi \to B$, with the total space $\Xi$ compact Hausdorff
and the trivial action of $\O (1)$ on the base $B$, which is locally  
$\O (1)$-equivariantly fibre homotopy equivalent
to a trivial bundle with fibre $S(\Rr^n)$ equipped with the antipodal
involution. 
We assume that the $\O (1)$-space $\Xi$ admits an $\O (1)$-equivariant
embedding as a subspace of some finite dimensional real $\O (1)$-module.

From a real line bundle $\lambda$ over $E$ we can form a spherical fibration
$\Xi_\lambda \to E$:
$$
\Xi_\lambda = (S(\lambda )\times_B \Xi )/\O (1),
$$
where $\O (1)$ acts as $\pm 1$ on $S(\lambda )$.

We consider an $\O (1)$-map $f : S(\lambda ) \to C_B(\Xi )$ to the fibrewise
cone on $\Xi$ and write $\tilde Z =\{ v\in S(\lambda )\st f(v)=0\}$
(where $0$ is the vertex of the cone in the fibre).
It determines a section $s$ of $C_E(\Xi_\lambda )$ with zero-set $Z$.
(If $\Xi=S(\xi )$ with the antipodal involution, then $C_B(\Xi )$ is the disc
bundle $D(\xi )$, $\Xi_\lambda = S(\lambda\otimes\xi )$ and
$C_E(\Xi_\lambda ) = D(\lambda\otimes\xi )$.)

The mod $2$ Euler class $e(\Xi_\lambda )\in H^n(E;\,\Ff_2)$
can be written as
$$
e(\Xi_\lambda ) = \sum_{i=0}^n e(\lambda )^i w_{n-i}(\Xi),
$$
where the classes $w_j(\Xi)\in H^j(B;\,\Ff_2)$ are defined by the universal
example in which $E=B\times P(\Rr^N)$ and $\lambda$ is the pullback of
the Hopf line bundle over $P(\Rr^N)$ for $N$ sufficiently large.

Then Proposition \ref{borsuk} generalizes as follows.
\Prop{Suppose that $a\in H^i(E;\,\Ff_2 )$ is a cohomology class 
such that $a\cdot e(\Xi_\lambda )\not=0\in H^{i+n}(E;\, \Ff_2)$.
Then $a$ restricts
to  a non-zero class in $H^i(Z ;\,\Ff_2 )$. 
}
\begin{proof}
This follows from the generalization,
included as Proposition \ref{loc} in the Appendix,
of \cite[Proposition 2.7]{borsuk}
in which $S(\xi )$ is replaced by an arbitrary spherical fibration.
\end{proof}

In the complex case we replace $\O (1)$ by $\U (1)$.
The spherical fibration $\Xi_\lambda$ is oriented and we have classes
$c_j(\Xi )\in H^{2j}(B;\,\Zz )$.

\Sect{Appendix: Euler classes for spherical fibrations}
We follow closely the account in \cite[Section 2]{borsuk}
to construct Euler classes with integral coefficients.
Corresponding definitions may be made, and theorems proved, for
$\Ff_p$-coefficients by purely notational changes.

Let $X$ be a compact ENR and let $\Xi \to X$ be a fibrewise space, with
$\Xi$ compact Hausdorff, which is locally fibre homotopy equivalent to
a trivial sphere bundle with fibre $S^{n-1}=S(\Rr^n)$.
From the spherical fibration $\Xi\to X$ we form the fibrewise cone
$C_X(\Xi )\to X$, which contains $\Xi$ as a subspace. We shall refer to the
vertex of the cone in any fibre as zero, written as `$0$'.
The non-zero elements of $C_X(\Xi )$ will be written as $[t,v]$,
where $t\in (0,1]$ and $v\in\Xi$.

We begin with the definition of the Euler class 
$e(\Xi )$ in the cohomology group $H^n(X;\,\Zz (\Xi ))$
with integral coefficients $\Zz (\Xi )$ 
twisted by the orientation bundle of $\Xi$.
There is a unique class
$u\in H^n(C_X(\Xi ),\Xi ;\, \Zz (\Xi ))$ restricting to the canonical generator
of each fibre $H^n(C (\Xi_x), \Xi_x;\, \Zz (\Xi_x))=\Zz$.
We define $e(\Xi )$ to be the pullback $z^*u$ by the
zero-section $z : (X,\emptyset )\to (C_X(\Xi) ,\Xi )$.

Let $s: X\to C_X(\Xi )$ be a section. We use the notation
$$
\Zero (s) =\{ x\in X \st s(x)=0\}
$$
for its {\it zero-set}.
It is a closed, so compact, subspace of $X$.

\Def{\label{csEuler}
Suppose that
$U\subseteq X$ is an open subset and that $s$ is a section of $C_X(\Xi )$ 
such that $\Zero (s)\cap U$ is compact.
Choose an open neighbourhood $W$ of $\Zero (s)\cap U$ such that $\cl{W}$ is
compact and contained in $U$.
There is an $\epsilon \in (0,1)$ such that for each
$x\in \cl{W}-W$ we have $s(x)=[t,v]$,
where $t\in [\epsilon, 1]$, $v\in\Xi_x$. 
For $x\in\cl{W}$, put 
$$
s'(x)= 
\begin{cases}
0&\text{if $s(x)=0$,}\\
[t/\epsilon ,v]&\text{if $s(x)=[t,v]$ with $t\in (0,\epsilon ]$, $v\in\Xi_x$,}
\\
[1 ,v]=v&\text{if $s(x)=[t,v]$ with $t\in [\epsilon,1]$, $v\in\Xi_x$,}
\end{cases}
$$
so that $s'$ gives a map $(\cl{W},\cl{W}-W)\to (C_X(\Xi ) ,\Xi)$.
The pullback 
of the Thom class $u\in H^n(C_X(\Xi ), \Xi ;\, \Zz (\Xi ))$
by this map $s'$ 
gives a class in the group
$H^n(\cl{W},\cl{W}-W;\, \Zz (\Xi ))$ which we may identify
with the cohomology $H^n_c(W;\,\Zz (\Xi ))$ of $W$ with compact supports.

The image of this class in $H^n_c(U;\,\Zz (\Xi))$ under the homomorphism
induced by the inclusion $W\into U$ is independent of the
choices made and will be called the
{\it Euler class with compact supports}
$$
e(s\, |\, U)\in H_c^n(U;\,\Zz (\Xi ))
$$
of the section $s$.
}
\Lem{\label{csprop}
{\rm (Properties of the Euler class with compact supports).}
\par\noindent
{\rm (i)}\quad
Suppose that $U'\subseteq U$ is an open subspace of $U$ such that
$\Zero (s)\cap U \subseteq U'$. Then $e(s\, |\, U) =i_! e(s\, |\, U')$,
where $i_! : H^n_c(U';\, \Zz (\Xi )) \to H^n_c(U;\, \Zz (\Xi ))$ is
induced by the inclusion $i :U'\into U$. 
\par\noindent
{\rm (ii)}\quad
Suppose that $U=U_1\sqcup U_2$ is a disjoint union of two open subspaces
$U_1$ and $U_2$. Then $e(s\, |\, U) = (i_1)_!e(s\, |\, U_1) +
(i_2)_!e(s\, |\, U_2)$, where $i_1: U_1\into U$ and $i_2: U_2\into U$
are the inclusion maps.
}
We shall use the Euler class with compact supports to localize the Euler class
of a spherical fibration to an arbitrarily small
neighbourhood of the zero-set of a section.
The basic result follows directly from the definition (and is, indeed,
a special case of the property (i) above).
\Lem{\label{basic}
Let $s$ be a section of $C_X(\Xi )$ with zero-set $\Zero (s)\subseteq U$,
where $U\subseteq X$ is open.
Then the Euler class $e(s\, |\, U)\in H^n_c(U;\,\Zz (\Xi ))$ with
compact supports maps under the homomorphism
$j_!$ induced by the inclusion $j:U\into X$
to $e(\Xi )\in H^n(X;\,\Zz (\Xi ))$.
}
This yields the basic result on the cohomology of the zero-set.
\Prop{\label{loc}
Let $s$ be a section of $C_X(\Xi )$ with zero-set $\Zero (s)$.
Suppose that $a\in H^i(X;\,\Zz )$ is a cohomology class that restricts
to $0\in H^i(\Zero (s) ;\,\Zz )$. Then 
$a\cdot e(\Xi )=0\in H^{i+n}(X;\, \Zz (\Xi ))$.
}
\begin{proof}
Choose an open neighbourhood $U$ of $\Zero (s)$ in $X$ such that
$a$ restricts to zero in $H^i(\cl{U};\,\Zz)$.
We have an Euler class with compact supports 
$e(s\, |\, U)\in H^n_c(U;\, \Zz (\Xi ))$
mapping to $e(\Xi )\in H^n(X;\, \Zz (\Xi ))$.
So $a\cdot e(s\, |\, U)$ maps to $a\cdot e(\Xi )$.
But $a$ restricts to $0$ in $H^i(\cl{U};\,\Zz )$ and hence
$a\cdot e(s\, |\, U)=0$.
\end{proof}
\Cor{\label{loccor}
Suppose that the kernel of multiplication by $e(\Xi )$:
$$
H^i(X;\,\Zz )\to H^{i+n}(X;\,\Zz (\Xi ))
$$
is zero. Then the restriction map
$$
H^i(X;\,\Zz ) \to H^i(\Zero (s);\, \Zz )
$$
is injective.
\qed
}
\Cor{\label{cdim}
Suppose that $s$ is a section of $C_X(\Xi )$ with zero-set
$\Zero (s)$. If there is a class $x\in H^i(X;\,\Zz )$ such that
$a\cdot e(\xi )\not=0$, then the covering dimension
of $\Zero (s)$ is greater than or equal to $i$. \qed
}
\Prop{\label{cover}
Let $X$ be a closed smooth manifold and let 
$Y\subseteq X$ be a closed submanifold of codimension $d$.
Suppose that $s$ is a section of $C_X(\Xi )$. If there is a class
$a\in H^i(Y;\,\Zz )$ such that $a\cdot e(\Xi |Y) \not=0\in H^{i+n}(Y;\,
\Zz (\Xi ))$, then the covering dimension of $\Zero (s)\subseteq X$
is at least $d+i$.
}
\begin{proof}
Writing $\nu$ for the normal bundle of $Y$ in $X$, choose a 
tubular neighbourhood $D(\nu )\into X$.
Now we have a map
$$
f: X' = S(\Rr\oplus\nu ) \to D(\nu )\into X,
$$
given by the projection $(t,v)\mapsto v$.
The section $s$ lifts to a section $s'$ of
$C_{X'}(f^*\Xi )$. And its zero-set $\Zero (s')$ is the pullback of
$\Zero (s) \cap D(\nu )$.
Let $u\in H^d(S(\Rr\oplus\nu );\,\Zz (\nu ))$ be the Thom class
of $\nu$.
Then $a'=u\cdot f^*a\in H^{d+i}(X';\, \Zz (\nu ))$ is non-zero.
Applying Corollary \ref{cdim} (with twisted coefficients) to
$s'$, we see that the covering dimension of
$\Zero (s')$ is at least $d+i$.

But $\Zero (s')$ is the union of two closed subspaces homeomorphic
to $\Zero (s)\cap D(\nu)$. It follows that $\Zero (s)\cap D(\nu )$
and {\it a fortiori} $\Zero (s)$ have covering dimension greater
than or equal to $d+i$.
\end{proof}
\Rem{If $n$ is large, then $\Xi$ is canonically equivalent to a free
$\O (1)$-spherical fibration.
See \cite[Theorem 3.10]{z2}.
For this canonical structure, the classes $w_j(\Xi )$ defined
in Section 5 satisfy $w_j(\Xi ) \cdot u = \Sq^ju$, where
$u\in H^n(C_X(\Xi ),\Xi ;\,\Ff_2)$ is the mod $2$ Thom class.
}

\medskip

\vfill\eject

\par\noindent{\bf Footnote.}
This paper originated in the observation by the first author that
two remarks (Remark 6.1 and Remark 6.2) in \cite{S1} were incompatible 
with a result of Stolz \cite{stephan} (see, also, \cite{z4})
on the level of real projective spaces.
The error lies in the proof of Corollaries 4.2 and 4.4, which
is modelled on the somewhat misleading approach in \cite{KM} to
the result stated there as Corollary 1.5.
We trust that the argument used in \cite[Corollary 1.5]{KM},
and also in
\cite[Theorem 1.3]{MPSS},
\cite[Theorem 1.5]{MPS} and \cite[Theorem 1.3]{S2},
will now be superseded by our proof of those theorems as
applications (Examples \ref{exs}) of Corollary \ref{simple}.

\end{document}